\documentclass[fleqn]{mat01}
\usepackage{times,mathtimy,amssymb,latexsym,amsbsy}
\begin{document}

\setcounter{page}{509} \firstpage{509}

\newtheorem{theore}{Theorem}
\renewcommand\thetheore{\arabic{section}.\arabic{theore}}
\newtheorem{theor}[theore]{\bf Theorem}
\newtheorem{rem}[theore]{Remark}
\newtheorem{propo}[theore]{\rm PROPOSITION}
\newtheorem{lem}[theore]{Lemma}
\newtheorem{definit}[theore]{\rm DEFINITION}
\newtheorem{coro}[theore]{\rm COROLLARY}
\newtheorem{exampl}[theore]{Example}
\newtheorem{case}{Case}

\def\corol{\trivlist \item[\hskip \labelsep{COROLLARY.}]}
\def\noteproof{\trivlist \item[\hskip \labelsep{\it Note added in Proof.}]}
\def\proof{\trivlist \item[\hskip \labelsep{\it Proof.}]}
\def\rema{\trivlist \item[\hskip \labelsep{\it Remark.}]}

\def\d{\mbox{\rm d}}
\def\e{\mbox{\rm e}}

\renewcommand{\theequation}{\thesection\arabic{equation}}

\title{Large time behaviour of solutions of a system of generalized Burgers
equation}

\markboth{K T~Joseph}{Asymptotic behaviour of solutions}

\author{K~T~JOSEPH}

\address{School of Mathematics, Tata Institute of Fundamental Research,
Homi Bhabha Road, Mumbai~400~005, India\\
\noindent E-mail: ktj@math.tifr.res.in}

\volume{115}

\mon{November}

\parts{4}

\pubyear{2005}

\Date{MS received 28 October 2004; revised 25 July 2005}

\begin{abstract}
In this paper we study the asymptotic behaviour of solutions of
a system of $N$ partial differential equations. When
$N = 1$ the equation reduces to the Burgers equation and was studied
by Hopf. We consider both the inviscid and viscous case and show
a new feature in the asymptotic behaviour.
\end{abstract}


\maketitle

\section{Introduction}

We consider the system of generalized Burgers equations for $N$
unknown variables $u = (u_1, u_2, \dots, u_N)$,
\begin{equation}
(u_j)_t + \sigma(c,u) (u_j)_{x}  = \frac{\nu}{2} (u_j)_{xx}, \quad
j = 1,2,\dots ,N,
\end{equation}
where $c = (c_1, c_2, \dots, c_N)$ is a constant vector in $R^N$
and $\sigma(c,u) = \sum_{k=1}^N c_k u_k$ is the usual inner
product in $R^N$. We study the solution of (1.1) with initial
conditions
\begin{equation}
u_j(x,0)=u_{0j}(x), \quad j = 1,2,\dots, N.
\end{equation}

When $\nu >0$, (1.1) is a system of nonlinear parabolic equation
describing the interplay between nonlinearity and diffusion, $\nu$
being the viscosity parameter.

When $\nu = 0$, the system (1.1) is hyperbolic with coinciding wave
speeds $\sigma(c, u)$, and the nonlinearity and the nonconservative form
makes the initial value problem complex. Indeed the system being
nonlinear, solution cannot be continued as a smooth solution even when
the initial data is smooth. Further if $N>1$, the product $\sigma
(c,u) \cdot (u_j)_{x}$ is nonconservative and does not make sense in the usual
distributional sense. The solution should be understood in a generalized
sense.

For the case $\nu > 0$, Joseph \cite{ktj02} used a generalized Hopf--Cole
transformation to linearize the system of equations $(1.1)$ and solve
(explicitly) with the conditions $(1.2)$ in terms of a family of
probability measures $\d\mu_{(x,t)}^{\nu}(y)$. These measures depend
on the initial data (1.2) in a nonlinear and nonlocal manner and takes the
form
\begin{equation}
\d\mu_{(x,t)}^{\nu}(y) = \frac{\e^{-\frac{1}{\nu}
\left[I(y)+\frac{(x-y)^2}{2t}\right]}\d y}{\int_{-\infty}^{\infty}
\e^{-\frac{1}{\nu} \left[I(y) + \frac{(x-y)^2}{2t}\right]}\d y},
\end{equation}
where
\begin{equation}
I(y) = \int_{0}^{y} \sigma(c,u_{0}(z)) \d z.
\end{equation}
When $\nu > 0$, the  solution of (1.1) and (1.2) was shown to be
\begin{equation}
u_{j}^\nu(x,t) = \int_{R^1}u_{0j}(y) \hbox{d}\mu_{(x,t)}^{\nu}(y),
\quad j = 1, 2, 3,\dots, N.
\end{equation}

Let $u_{j}^{\nu}(x,t)$ be the solution of (1.1) and (1.2) given by
(1.5). It was shown in \cite{ktj03} that when $u_{0j}$ is Lipschitz continuous,
for each $t>0$, except for a countable $x$ the limits
\begin{equation*}
u_j(x,t) = \lim_{\nu \rightarrow 0} u_{j}^{\nu}(x,t)
\end{equation*}
exist and is given by the formula
\begin{equation}
u_j(x,t) = u_{0j}(y(x,t)), \quad j=1,2,\dots,N
\end{equation}
where $y(x,t)$ is a minimizer of
\begin{equation}
\min_{-\infty <y<\infty} (I(y) + (x-y)^2/2t)
\end{equation}
and $I(x)$ is given by (1.4).

Using some ideas from the earlier works of Joseph \cite{ktj93} and
LeFloch \cite{Le90}, $u_j(x,t)_{j = 1,2,\dots,N}$ was shown to be
a solution of an inviscid case $(\nu = 0)$ in (1.1) with initial
data from (1.2), the nonconservative product was justified in the
sense of Volpert product \cite{Vol67} (see Dal Maso, LeFloch and
Murat \cite{Dal95} for a generalization of Volpert product). In
\cite{ktj02} solution for general initial data is constructed in
the sense of Colombeau \cite{Col84}.

The aim of the present note is to study the asymptotic behaviour
of the solution for the parabolic (viscous) case as well as the
hyperbolic (inviscid) case. Study of asymptotic behaviour of
solutions is important, see \cite{Sah87} and the references
therein for the parabolic case and \cite{Lax57} for the inviscid
case. When $N = 1$, this system is the celebrated Burgers equation
and explicit solution and its asymptotic behavior as $t$ tends to
infinity and diffusion parameter $\nu \rightarrow 0$ was studied
by Hopf \cite{H50}. We show that Hopf's analysis give asymptotic
form of the solution for the viscous case. When $N=1$ and $c_1
\neq 0$, it is well-known from the work of Lax \cite{Lax57} that
for solution of inviscid Burgers equation with initial data
supported in the compact interval $[-\ell,\ell]$, $\ell>0$ the
solution decays at the rate $O(t^{-\frac{1}{2}})$ and support
spreads at a rate $O(t^{\frac{1}{2}})$ for large time. From an
explicit computation we will show that the decay rate is not true
in general for the present case, but still the support spread at
the same rate. We start with the viscous case.

\section{Asymptotic behaviour with viscous term}

\setcounter{equation}{0}

In this section we study the asymptotic behavior of solution of
(1.1) and (1.2) when $\nu > 0$ and fixed. On the initial
conditions $u_{0j}(x), j=1,2,\dots,N$, assume that $\lim_{x
\rightarrow \infty} I(x) = I(\infty)$, $\lim_{x \rightarrow
-\infty} I(x) = I(-\infty)$, $\lim_{x \rightarrow \infty}
u_{0j}(x) = u_{0j}(\infty)$, $\lim_{x \rightarrow -\infty}
u_{0j}(x) = u_{0j}(-\infty)$ exists and is finite. With these
assumptions, we shall prove the following result.

\begin{theor}[\!]
The solution
$u^{\nu}(x,t)=(u_1^\nu(x,t),u_2^\nu(x,t),\dots,u_N^\nu(x,t))$ of
{\rm (1.1)} and {\rm (1.2)} has the following asymptotic behaviour
as $t$ tends to infinity\hbox{{\rm :}}
\begin{equation}
\hskip -3pc u_j^{\nu}(x,t) \approx
\frac{u_{0j}(\infty) \e^{\frac{-I(\infty)}{\nu}}
\int_{-\infty}^{x/\sqrt(t\nu)} \e^{-\frac{y^2}{2}}\d y +
u_{0j}(-\infty) \e^{\frac{-I(-\infty)}{\nu}}
 \int_{x/\sqrt(t\nu)}^{\infty} \e^{-\frac{y^2}{2}}\d y}{
\e^{\frac{-I(\infty)}{\nu}} \int_{-\infty}^{x/\sqrt(t\nu)}
\e^{-\frac{y^2}{2}}\d y + \e^{\frac{-I(-\infty)}{\nu}}
 \int_{x/\sqrt(t\nu)}^{\infty} \e^{-\frac{y^2}{2}}\d y}.
\end{equation}
\end{theor}

\begin{proof}
First we note that the solution $u_j^\nu(x,t),j=1,2,\dots,N$ of
(1.1) and (1.2) is given by (1.5) where the measure
$\d\mu_{(x,t)}^{\nu}(y)$ is given by eqs~(1.3) and (1.4). Writing
explicitly the formula, we get
\begin{equation*}
u_{j}^\nu(x,t)=\frac{\int_{-\infty}^{\infty} u_{0j}(y)
\e^{-\frac{1}{\nu} \big[I(y)+\frac{(x-y)^2}{2t}\big]}\d
y}{\int_{-\infty}^{\infty} \e^{-\frac{1}{\nu} \left[I(y) +
\frac{(x-y)^2}{2t}\right]}\d y}.
\end{equation*}

Setting $\xi = x/\sqrt{\nu  t}$, and then making a change of
variable $z = \frac{\sqrt{\nu t}\xi - y}{\sqrt{\nu t}}$ and
renaming $z$ as $y$, we get
\begin{equation}
u^{\nu}_{j}(x,t) = \frac{\int_{-\infty}^{\infty} u_{0j}(\sqrt{\nu
t}(\xi - y) \e^{-\left[\frac{I(\sqrt{\nu  t}(\xi - y)}{\nu} +
y^2/2\right]} \d y}{\int_{-\infty}^{\infty}
\e^{-\left[\frac{I(\sqrt{\nu  t}(\xi - y))}{\nu} + y^2/2\right]}
\d y}.
\end{equation}

Now we split the integral in (2.2) in the following manner:
\begin{align}
&\int_{-\infty}^{\infty} u_{0j}(\sqrt{\nu t}(\xi - y))
\e^{-\left[\frac{I(\sqrt{\nu t}(\xi - y))}{\nu} +y^2/2\right]} \d y\nonumber\\[.5pc]
&\quad = \int_{-\infty}^{\xi -\delta} u_{0j}(\sqrt{\nu t}(\xi -
y))\e^{-\left[\frac{I(\sqrt{\nu t}(\xi - y))}{\nu}+y^2/2\right]} \d y\nonumber\\[.5pc]
&\qquad + \int_{\xi +\delta}^{\infty} u_{0j}(\sqrt{\nu t}(\xi -
y))\e^{-\left[\frac{I(\sqrt{\nu t}(\xi - y))}{\nu}+y^2/2\right]} \d y\nonumber\\[.5pc]
&\qquad +\int_{\xi-\delta}^{\xi+\delta} u_{0j}(\sqrt{\nu t}(\xi -
y)) \e^{-\left[\frac{I(\sqrt{\nu t}(\xi - y))}{\nu}+y^2/2\right]}
\d y.
\end{align}
Now we fix $\delta >0$ and study each of these integrals as $t$ tends to
infinity. We have under the assumptions of the theorem on $I(x)$ and
$u_{0j}(x)$, as $t$ tends to infinity:
\begin{align*}
\hskip -3pc &\int_{-\infty}^{\xi - \delta} u_{0j}(\sqrt{\nu t}(\xi
- y)) \e^{-\left[\frac{I(\sqrt{\nu t}(\xi - y))}{\nu}
+y^2/2\right]} \d y \approx
\e^{-\frac{I(+\infty)}{\nu}}u_{0j}(\infty)\int_{-\infty}^{\xi
-\delta} \e^{-y^2/2} \d y.\\[.5pc]
\hskip -3pc &\int_{\xi + \delta}^{\infty} u_{0j}(\sqrt{\nu t}(\xi
- y)) \e^{-\left[\frac{I(\sqrt{\nu t}(\xi - y))}{\nu}
+y^2/2\right]} \d y \approx
\e^{-\frac{I(-\infty)}{\nu}}u_{0j}(-\infty)\int_{\xi +
\delta}^{\infty} \e^{-y^2/2} \d y, \\[.5pc]
\hskip -3pc &\limsup_{t\rightarrow \infty} \left\vert\int_{\xi -
\delta}^{\xi + \delta} u_{0j}(\sqrt{\nu t}(\xi - y))
\e^{-\left[\frac{I(\sqrt{\nu t}(\xi - y))}{\nu} +y^2/2\right]} \d
y\right\vert = O(\delta).
\end{align*}
Now let $t$ tend to infinity and then $\delta$ tend to $0$
in (2.3). We get
\begin{align}
\hskip -3pc \begin{split} \int_{-\infty}^{\infty} u_{0j}(\sqrt{\nu
t}(\xi - y)) \e^{-\left[I(\sqrt{\nu t}(\xi - y)) \nu
+y^2/2\right]} \d y &\approx
\e^{-\frac{I(+\infty)}{\nu}}u_{0j}(\infty)\int_{-\infty}^{\xi}
\e^{-y^2/2} \d y\\[.5pc]
&\quad + \e^{-\frac{I(-\infty)}{\nu}}u_{0j}(-\infty)\int_{\xi}^{\infty}
\e^{-y^2/2} \d y.
\end{split}
\end{align}
Similarly
\begin{equation}
\hskip -2pc \int_{-\infty}^{\infty} \e^{-\left[\frac{I(\sqrt{\nu
t}(\xi - y))}{\nu}+y^2/2\right]} \d y \approx
\e^{\frac{-I(+\infty)}{\nu}} \int_{-\infty}^{\xi} \e^{-y^2/2} \d y
+ \e^{\frac{-I(-\infty)}{\nu}} \int_{\xi}^{\infty} \e^{-y^2/2} \d
y.
\end{equation}
We observe that due to our assumption on $I(x)$, this limit in (2.5) is a
postive real number and hence letting $t$ tend to infinity in (2.2) and
using (2.4) and (2.5) we get the result (2.1). The proof of the theorem
is complete. \hfill $\Box$
\end{proof}

\begin{rema}
{\rm An interesting case here is when the initial data $u_{0j}, j
= 1,2,\dots,N$ satisfies the following conditions.
$u_{0j}(\infty)$, $u_{0j}(-\infty)$ are nonzero and there is a
cancellation in $\sigma(c,u_0)(x)=\sum_1^N c_k u_{0k}(x)$ so that
this quantity is integrable.}
\end{rema}

\section{Asymptotic behaviour of solutions of generalized Hopf equation}

\setcounter{equation}{0}
\setcounter{theore}{0}

In this section we study the asymptotic behaviour of vanishing viscosity
solutions of
\begin{equation}
(u_j)_t + \sigma(c,u)(u_j)_x = 0
\end{equation}
with initial data
\begin{equation}
u_j(x,0) = u_{0j}(x),
\end{equation}
for $j=1,2,\dots, N$. We recall the definition of $\sigma(c,u)$
namely $\sigma(c,u)=\sum_1^N c_k u_k$ where
$c=(c_1,c_2,\dots,c_N)$ a given constant vector and $u =
(u_1,u_2,\dots,u_N)$, the unknown vector variable.

It is easy to see that $u_j^\nu(x,t),j = 1,2,\dots,N$ is a
solution of (1.1) and (1.2) iff $\sigma^\nu(x,t):=
\sigma(c,u^\nu(x,t))$ satisfies the Burgers equation
\begin{equation*}
\sigma_t + \left(\frac{\sigma^2}{2} \right)_x = \frac{\nu}{2}
\sigma_{xx}
\end{equation*}
with initial condition
\begin{equation*}
\sigma(x,0) = \sigma(c,u_0)(x).
\end{equation*}
By Hopf--Cole transformation, the solution can be written in the form
\begin{equation*}
\sigma^\nu(x,t)=\int_{R^1}\sigma(c,u_0(y)) \d\mu_{(x,t)}^{\nu}(y).
\end{equation*}
As in \cite{ktj03} it is easy to see that, when $u_{0j}$ is Lipschitz continuous,
for each $t > 0$, except for a countable number of points of $x$, the limit
\begin{equation*}
\sigma(x,t) = \lim_{\nu \rightarrow 0} \sigma^{\nu}(x,t)
\end{equation*}
exists and is given by the formula
\begin{equation*}
\sigma(x,t)= \sigma(c,u_0(y(x,t)))= \sigma(c,u),
\end{equation*}
where $y(x,t)$ is a minimizer of
\begin{equation*}
\min_{-\infty <y<\infty} (I(y) + (x-y)^2/2t)
\end{equation*}
and $I(x)$ is given by (1.4). Here we remark that this formula is
slightly different from that of Hopf \cite{H50}. Once
$\sigma(c,u)$ is known, (3.1) can be treated as $N$ scalar linear
equation with discontinuous coefficient. For any function $a(x,t)$
on $[(x,t)\hbox{\rm :}\ -\infty < x < \infty, t>0]$,
$\|a(.,t)\|_\infty$ denotes the essential supremum of $a(x,t)$
with respect to the space variable $x$ keeping the time variable
$t \geq 0$ fixed. We define the left boundary curve $x=s^-(t)$ and
right boundary curve $x = s^+(t)$ of the support of $a(x,t)$ where
$s^-(t) = \sup[y\hbox{\rm :}\ a(x,t) = 0$ for all $x < y]$ and
$s^+(t) = \inf[y\hbox{\rm :}\ a(x,t)=0$ for all $x > y]$. It is
well-known from \cite{Lax57} that when the initial data is
supported in $[-l,l]$, being solution of the inviscid Burgers
equation, $\sigma(c,u)$ have the following estimates. There exist
constants $A>0$ and $C>0$ which depend on $l$ and
$\|u_0\|_{\infty}$, such that
\begin{align}
-l - A t^{\frac{1}{2}} \leq s^{-}(t) &\leq s^{+}(t) \leq l + A
t^{\frac{1}{2}},\nonumber \\
\|\sigma(.,t)\|_{\infty} &\leq C t^{\frac{-1}{2}}.
\end{align}
We shall prove that for $u_j, j=1,2,\dots,N$ again the same
estimate holds for the spread of support but the decay result is
not valid in general. First we have the following theorem.

\begin{theor}[\!]
Let $u_j(x,t), j = 1,2,\dots,N$ be the vanishing viscosity
solution of {\rm (3.1)} and {\rm (3.2)} with initial data
supported in $[-l,l]$. Let $x=s_j^-(t)$ and $x=s_j^+(t)$ are the
support curves for $u_j(x,t)$. Then there exists a constant $A>0$
which depend on $l$ and $\|u_0\|_\infty${\rm ,} so that the
following estimate holds for $t \gg 1${\rm ,}
\begin{equation*}
-l - A t^{\frac{1}{2}} \leq s_j^{-}(t) \leq s_j^{+}(t) \leq l + A
t^{\frac{1}{2}}.
\end{equation*}
\end{theor}

\begin{proof}
The proof easily follows from the fact that the characteristic
speed $\sigma(c,u) = 0$ outside the region $ -l - A
t^{\frac{1}{2}} \leq x \leq l + A t^{\frac{1}{2}}$. So in this
region the characteristics connecting the points $(x,t)$ to a base
point $(y,0)$ are parallel to the $t$-axis and hence is of the
form $x = y$ and since the solution is constant along the
characteristics, we have $u(x,t) = u_0(y) = u_0(x)$. When $(x,t)$
is outside $ -l - A t^{\frac{1}{2}} \leq x \leq l + A
t^{\frac{1}{2}}${\rm ,} the $x$ co-ordinate of the base point of
the characteristic{\rm ,} $y${\rm ,} lies outside $[-l,l]$ where
$u_0(y)$ is zero and hence $u(x,t)$ is zero. The proof of the
theorem is complete. \hfill $\Box$
\end{proof}

Next we show that solution does not decay in general, by giving an
example. Here we construct vanishing viscosity solution of (3.1) with
initial data of the special form which is supported in a compact set, namely
\begin{equation}
u_j(x,0) = \begin{cases}
0, &{\rm if} \ x < - l\\
u_{0j}, & {\rm if} \ -l < x < l\\
0, &{\rm if} \ x >l
\end{cases}
\end{equation}
where $u_{0j}$ is a constant and $l$ is a positive real number.
Let $u_j^\nu(x,t)_{j=1,2,\dots,N}$ be the solution of (1.1) with
the initial data (3.4). We shall prove the following.

\begin{theor}[\!]
Let $\sigma_0=\sum_{k=1}^N c_j u_{0j}${\rm ,} then $u_j(x,t)=\lim_{\nu
\rightarrow 0} u_j^\nu(x,t)$ exists and takes the
following form\hbox{{\rm :}}

When $\sigma_0 = 0${\rm ,}
\begin{equation}
u_j(x,t) = \begin{cases} 0, &{\rm if} \ x <-l\\
u_{0j}, &{\rm if} \ -l <x < l\\
0, & {\rm if} \ x >l
\end{cases}.
\end{equation}

When $\sigma_0 < 0${\rm ,}
\begin{equation}
u_j(x,t) = \begin{cases}
0, &{\rm if} \ x < \frac{\sigma_0}{2}\cdot t - l, t< \frac{-4l}{\sigma_0}\\[.3pc]
u_{0j}, &{\rm if} \ \frac{\sigma_0}{2}t - l < x < \sigma_0 t + l,
t<\frac{-4 l}{\sigma_0}\\[.3pc]
\frac{u_{jl}}{\sigma_0}\cdot \frac{x-l}{t}, &{\rm if} \ \sigma_0 t
< x < l,
t<\frac{-4 l}{\sigma_0}\\[.3pc]
0, &{\rm if} \ x <l - (-4 l \sigma_0 t)^{\frac{1}{2}}, t >
\frac{-4l}{\sigma_0}\\[.3pc]
\frac{u_{0j}}{\sigma_0}\cdot \frac{x-l}{t}, &{\rm if} \
l - (-4 l \sigma_0 t)^{\frac{1}{2}} < x < l, t > \frac{-4 l}{\sigma_0}\\[.3pc]
0, &{\rm if} \ x > l
\end{cases}.
\end{equation}

When $\sigma_0 >0${\rm ,}
\begin{equation}
\lim_{\nu \rightarrow 0}u_j^\nu(x,t) = \begin{cases}
0, &{\rm if} \ x > \frac{\sigma_0}{2}\cdot t + l , t< \frac{4l}{\sigma_0}\\[.3pc]
u_{0j}, &{\rm if} \ \sigma_0 t - l <x < \frac{\sigma_0}{2}t + l, t<\frac{4
l}{\sigma_0}\\[.3pc]
\frac{u_{0j}}{\sigma_0}\cdot \frac{x+l}{t}, & {\rm if}\ -l < x <
\sigma_0 t
-l, t<\frac{4 l}{\sigma_0}\\[.3pc]
0, &{\rm if} \ x > l +(4 l \sigma_0 t)^{\frac{1}{2}}, t > \frac{4l}{\sigma_0}\\[.3pc]
\frac{u_{0j}}{\sigma_0}\cdot \frac{x+l}{t}, &{\rm if} \ -l < x < l
+ (4 l \sigma_0 t)^{\frac{1}{2}}, t > \frac{4 l}{\sigma_0}\\[.3pc]
0, &{\rm if}\  x <-l
\end{cases}.
\end{equation}
\end{theor}

\begin{proof}
To prove (3.5)--(3.7) we use the formula (1.5) to get explicit solution
of (1.1) and (3.4) in the form
\begin{equation}
u_j^\nu(x,t) = \frac{u_{j0} \int_{-l}^l
\e^{-\frac{1}{\nu}\big[\frac{(x-y)^2}{2t}+ \sigma_0 y\big]}\d
y}{\int_{-\infty}^{-l} \e^{-\frac{(x-y)^2}{2t\nu}}\d y +
\int_{-l}^l \e^{-\frac{1}{\nu}\big[\frac{(x-y)^2}{2t}+\sigma_0
y\big]}\d y + \int_{l}^{\infty} \e^{-\frac{(x-y)^2}{2t\nu}}\d y}
\end{equation}
which can be written in the form
\begin{equation}
u_j^\nu(x,t) = \frac{u_{j0} \int_{-l}^l
\e^{-\frac{1}{\nu}\big[\frac{(x-y)^2}{2t}+ \sigma_0 y\big]}\d
y}{(2t\pi \nu)^{\frac{1}{2}} + \int_{-l}^l
\e^{-\frac{1}{\nu}\big[\frac{(x-y)^2}{2t}+\sigma_0 y\big]}\d y +
\int_{-l}^{l} \e^{-\frac{(x-y)^2}{2t\nu}}\d y}.
\end{equation}
To study the limit, we rewrite this formula in a convenient form
by introducing the functions
\begin{equation}
A_{l,\sigma_0}^\nu(x,t)=(2t\nu)^{\frac{1}{2}} \e^{\frac{\sigma_0^2
t}{2\nu} -\frac{\sigma_0 x}{\nu}}
erfc\left(\frac{t\sigma_0-x-l}{(2t\nu)^{\frac{1}{2}}}\right)
\end{equation}
and
\begin{equation}
B_{l,\sigma_0}^\nu(x,t)=(2t\nu)^{\frac{1}{2}} \e^{\frac{\sigma_0^2 t}{2\nu}
-\frac{\sigma_0 x}{\nu}} erfc\left(\frac{t\sigma_0-x+l}{(2t\nu)^{\frac{1}{2}}}\right),
\end{equation}
where
\begin{equation}
erfc(y)=\int_y^\infty \e^{-y^2} \d y.
\end{equation}
We can rewrite (3.9) as
\begin{equation}
u_j^\nu(x,t) = \frac{u_{0j}(A_{l,\sigma_0}^\nu(x,t)-
B_{l,\sigma_0}^\nu(x,t))}{ (2 t \pi
\nu)^{\frac{1}{2}} + A_{l,\sigma_0}^\nu(x,t)-B_{l,\sigma_0}^\nu(x,t)+
A_{l,0}^\nu(x,t)-B_{l,0}^\nu(x,t)}.
\end{equation}
Using the asymptotic expansions of the {\it erfc}, namely,
\begin{equation*}
{\it erfc(y)} = \left(\frac{1}{2y}-\frac{1}{4y^3}
+ o\left(\frac{1}{y^3}\right)\right) \e^{-y^2}, \quad y \rightarrow \infty
\end{equation*}
and
\begin{equation*}
{\it erfc(-y)} =
(\pi)^{\frac{1}{2}} - \left(\frac{1}{2y}-\frac{1}{4y^3} + o \left(\frac{1}{y^3}\right)\right)
\e^{-y^2},  \quad y \rightarrow \infty
\end{equation*}
in (3.10) and (3.11) we have the following as $\nu\rightarrow 0$:
\begin{equation}
\hskip -4pc A_{l,\sigma_0}^{\nu}(x,t) \approx \begin{cases}
\frac{(t\nu)}{-l-x+\sigma_0t} \e^{-\frac{x^2}{2\nu t}}, &{\rm if} \ -l-x+\sigma_0t>0\\
(\frac{\pi t \nu}{2})^{\frac{1}{2}} \e^{\frac{\sigma_0^2
t}{2\nu}-\frac{\sigma_0x}{\nu}}, &{\rm if} \ -l-x+\sigma_0t=0\\
(2\pi t \nu)^{\frac{1}{2}} \e^{\frac{\sigma_0^2
t}{2\nu}-\frac{\sigma_0x}{\nu}} + \frac{(t\nu)}{-l-x+\sigma_0t}
\e^{-\frac{x^2}{2\nu t}}, &{\rm if} \
-l-x+\sigma_0 t < 0 \end{cases}.
\end{equation}\vspace{-1pc}
\begin{equation}
\hskip -4pc B_{l,\sigma_0}^{\nu}(x,t) \approx
\begin{cases}\frac{(t\nu)}{l-x+\sigma_0t}
\e^{-\frac{x^2}{2\nu t}}, &{\rm if} \ l-x+\sigma_0t > 0\\
(\frac{\pi t \nu}{2})^{\frac{1}{2}} \e^{\frac{\sigma_0^2
t}{2\nu}-\frac{\sigma_0x}{\nu}}, &{\rm if} \ l-x+\sigma_0t = 0\\
(2\pi t \nu)^{\frac{1}{2}} \e^{\frac{\sigma_0^2
t}{2\nu}-\frac{\sigma_0x}{\nu}} +
\frac{(t\nu)}{l-x+\sigma_0t} \e^{-\frac{x^2}{2\nu t}}, &{\rm if} \
l-x+\sigma_0 t < 0
\end{cases}.
\end{equation}

It is straightforward to check the formulas (3.5)--(3.7) using
(3.14) and (3.15) in (3.13).

When $\sigma_0 =0$, (3.13) becomes
\begin{equation}
u_j^\nu(x,t) = \frac{u_{0j}(A_{l,0}^\nu(x,t)-
B_{l,0}^\nu(x,t))}{ (2 t \pi \nu)^{\frac{1}{2}}}.
\end{equation}
Now take the region $x<-l$. Then $-l-x>0$ and $l-x>0$, and using
(3.14) and (3.15) in (3.16) we get
\begin{equation*}
u_j^\nu(x,t) \approx u_{0j}
\frac{\frac{t\nu}{-l-x} \e^{\frac{-x^2}{2t\nu}}-\frac{t\nu}{l-x} \e^{\frac{-x^2}{2t\nu}}}
{(2\pi t \nu)^{\frac{1}{2}}}
\end{equation*}
and hence we have
\begin{equation*}
\lim_{\nu\rightarrow 0}u_j^\nu(x,t)=0, \quad x<-l.
\end{equation*}
In the region $-l<x<l$, $-l-x<0, l-x>0$,
\begin{equation*}
u_j^\nu(x,t) \approx u_{0j}
\frac{(2\pi t \nu)^{\frac{1}{2}}+ \frac{t\nu}{-l-x} \e^{\frac{-x^2}{2t\nu}}-\frac{t\nu}{l-x}
\e^{\frac{-x^2}{2t\nu}}}{(2\pi t \nu)^{\frac{1}{2}}}.
\end{equation*}
Thus we get
\begin{equation*}
\lim_{\nu\rightarrow 0}u_j^\nu(x,t)=u_{0j}, \quad -l<x<l.
\end{equation*}
In the region $x>l$, we have $-l-x<0,l-x<0$ and
\begin{equation*}
u_j^\nu(x,t) \approx u_{0j}
\frac{(2\pi t\nu)^{\frac{1}{2}}+\frac{t\nu}{-l-x} \e^{\frac{-x^2}{2t\nu}}
-(2\pi t\nu)^{\frac{1}{2}}-\frac{t\nu}{l-x} \e^{\frac{-x^2}{2t\nu}}}{(2\pi t \nu)^{\frac{1}{2}}}
\end{equation*}
and hence it follows that
\begin{equation*}
\lim_{\nu\rightarrow 0}u_j^\nu(x,t)=0, \quad x > l.
\end{equation*}
This completes the proof of the theorem for $\sigma_0 =0$. The case
$\sigma\neq 0$ is similar and is omitted. \hfill $\Box$
\end{proof}

\begin{rema}
{\rm Thus if we take the initial data (3.4) which has compact
support the decay of the vanishing viscosity solution $u_j(x,t),
j=1,2,\dots, N$ depends on the initial speed $\sigma_0 = \sum_1^N
c_ku_{0k}$. If $\sigma_0 =0$, the vanishing viscosity solution
does not decay. Indeed
\begin{equation*}
\sup_{x\in R^1}|u_j(x,t)| = |u_0j|.
\end{equation*}
On the other hand, from the above theorem it follows that for the case
$\sigma_0\neq 0$, the solution decays, namely
\begin{equation*}
\sup_{x\in R^1}|u_j(x,t)| = O\left(t^{-{\frac{1}{2}}}\right).
\end{equation*}

We conclude with a remark on the solution for (3.1) with
the Riemann type initial data
\begin{equation}
(u_j)(x,0) = \begin{cases} u_{jL}, &{\rm if} \ x < 0\\
u_{jR}, &{\rm if} \ x > 0,\end{cases}
\end{equation}
where $u_{jL}$ and $u_{jR}$ are constants for $j=1,2,\dots ,n$.

Let $\sigma_L=\sum_{k=1}^{n}c_k u_{jL}$ and $\sigma_R=\sum_{k=1}^{n}c_k
u_{kR}$. Then the vanishing viscosity solution for the Riemann problem
(3.1) and (3.17) takes the following form \cite{ktj02}.

When $\sigma_L<\sigma_R$,
\begin{equation*}
u_j(x,t) = \begin{cases} u_{jL}, &{\rm if} \ x
\leq \sigma_Lt\\[.2pc]
\frac{u_{jR}-u_{jL}}{\sigma_R-\sigma_L}\cdot \frac{x}{t} +
\frac{u_{jL} \sigma_R - u_{jR} \sigma_L}{\sigma_R - \sigma_L},
&{\rm if} \ \sigma_{L}
t < x < \sigma_{R}t.\\[.2pc]
u_{jR}, &{\rm if} \ x\geq\sigma_Rt\end{cases}
\end{equation*}
When $\sigma_L=\sigma_R =\sigma$,
\begin{equation*}
u_j(x,t) = \begin{cases} u_{jL}, & {\rm if} \ x
< \sigma\cdot t\\[.2pc]
u_{jR}, &{\rm if} \ x > \sigma\cdot t.\end{cases}
\end{equation*}
When $\sigma_L>\sigma_R$,
\begin{equation*}
u_j(x,t) = \begin{cases} u_{jL}, &{\rm if} \ x
< \frac{\sigma_L +\sigma_R}{2}\cdot t\\[.3pc]
\frac{u_{jL}+u_{jR}}{2} &{\rm if} \ x
= \frac{\sigma_L +\sigma_R}{2}\cdot t\\[.3pc]
u_{jR}, &{\rm if} \ x>\frac{\sigma_L +\sigma_R}{2} \cdot t
\end{cases}.
\end{equation*}}
\end{rema}

\section*{Acknowledgement}

This work is supported by a grant (No.~2601--2) from the
Indo-French Centre for the Promotion of Advanced Research, IFCPAR
(Centre Franco-Indien pour la Promotion de la Recherche Avancee,
CEFIPRA), New Delhi.

\end{document}